\documentclass[12pt,leqno]{amsart}
\usepackage{amssymb,verbatim,enumerate,ifthen,cite}
\usepackage[mathscr]{eucal}
\usepackage[utf8]{inputenc}
\usepackage[T1]{fontenc}
\def\N{\mathbb{N}}
\def\R{\mathbb{R}}

\def\K{\mathbb{K}}

\def\C{\mathscr{C}}
\def\CC{\mathbb{C}}
\def\K{\mathbb{K}}

\def\cH{\mathscr{H}}

\long\def\comment#1{}

\newtheorem{theorem}{Theorem}[section]
\newtheorem*{theorem*}{Theorem}
\def\Thm#1#2{\ifthenelse{\equal{#1}{*}}{\begin{theorem*}#2\end{theorem*}}
             {\begin{theorem}\label{T#1}#2\end{theorem}}}

\def\thm#1{Theorem~\ref{T#1}}

\newtheorem{proposition}[theorem]{Proposition}
\newtheorem*{proposition*}{Proposition}
\def\Prp#1#2{\ifthenelse{\equal{#1}{*}}{\begin{proposition*}#2\end{proposition*}}
{\begin{proposition}\label{P#1}#2\end{proposition}}}

\newtheorem{corollary}[theorem]{Corollary}
\newtheorem*{corollary*}{Corollary}
\def\Cor#1#2{\ifthenelse{\equal{#1}{*}}{\begin{corollary*}#2\end{corollary*}}
             {\begin{corollary}\label{C#1}#2\end{corollary}}}
\def\cor#1{Corollary~\ref{C#1}}

\newtheorem{lemma}[theorem]{Lemma}
\newtheorem*{lemma*}{Lemma}
\def\Lem#1#2{\ifthenelse{\equal{#1}{*}}{\begin{lemma*}#2\end{lemma*}}
             {\begin{lemma}\label{L#1}#2\end{lemma}}}
\def\lem#1{Lemma~\ref{L#1}}

\theoremstyle{definition}
\newtheorem{remark}[theorem]{Remark}
\newtheorem*{remark*}{Remark}
\def\Rem#1#2{\ifthenelse{\equal{#1}{*}}{\begin{remark}\rm #2\end{remark}}
             {\begin{remark}\label{R#1}\rm #2\end{remark}}}

\newtheorem{example}[theorem]{Example}
\newtheorem*{example*}{Example}
\long\def\Exa#1#2{\ifthenelse{\equal{#1}{*}}{\begin{example*}\rm #2\end{example*}}
             {\begin{example}\label{Ex#1}\rm #2\end{example}}}
\def\exa#1{Example~\ref{Ex#1}}

\def\eq#1{{\rm(\ref{E#1})}}
\def\Eq#1#2{\ifthenelse{\equal{#1}{*}}
  {\begin{align*}#2\end{align*}}
  {\begin{equation}\begin{aligned}\label{E#1}#2\end{aligned}\end{equation}}}

\allowdisplaybreaks
\begin{document}
\begin{flushright}
\end{flushright}
\vspace{5mm}

\date{\today}

\title{Hermite-type interpolation in terms of exponential polynomials}

\author[A. H. Ali]{Ali Hasan Ali}
\address[A. H. Ali]{Doctoral School of Mathematical and Computational Sciences, University of Debrecen, H-4002 Debrecen, Pf.\ 400, Hungary}
\email{ali.hasan@science.unideb.hu}
\author[Zs. P\'ales]{Zsolt P\'ales}
\address[Zs. P\'ales]{Institute of Mathematics, University of Debrecen, 
H-4002 Debrecen, Pf.\ 400, Hungary}
\email{pales@science.unideb.hu}


\thanks{The research of the second author was supported by the K-134191 NKFIH Grant.}

\begin{abstract}
Motivated by classical results of approximation theory, we define an Hermite-type interpolation in terms of $n$-dimensional subspaces of the space of $n$ times continuously differentiable functions. In the main result of this paper, we establish an error term in integral form for this interpolation in the case when the $n$-dimensional subspace is the kernel of an $n$th order linear differential operator with constant coefficients. Several corollaries are deduced illustrating the applicability of this result.
\end{abstract}

\maketitle

\section{Introduction} 

In what follows, $I$ will always denote a nonempty open subinterval of the real line. The symbol $\C(I)$ will stand for the linear space of continuous complex-valued functions defined on $I$. If additionally $n\in\N$, then $\C^n(I)$ will denote the space of $n$-times continuously differentiable complex-valued functions defined on $I$.

To describe the classical result about Hermite interpolation, for $\ell,n\in\N$, consider the set
\Eq{*}{
  \cH_{\ell,n}(I):=\Big\{(\pmb{{a}},\pmb{n})\mid\, 
  & \pmb{{a}}=({a}_1,\dots,{a}_\ell)\in I^\ell \mbox{ with } {a}_1<\dots<{a}_\ell  \mbox{ and}\\
  & \pmb{n}=(n_1,\dots,n_\ell)\in\N^\ell \mbox{ with } n=n_1+\dots+n_\ell\Big\}.
}
Given an interpolational system $(\pmb{{a}},\pmb{n})=\big(({a}_1,\dots,{a}_\ell),(n_1,\dots,n_\ell)\big)\in\cH_{\ell,n}(I)$, the number $n$ will be called its \emph{dimension}. The order $\mathcal{O}(\pmb{{a}},\pmb{n})$ of $(\pmb{{a}},\pmb{n})$ is defined to be number $\max(n_1,\dots,n_\ell)-1$.

Let $(\pmb{{a}},\pmb{n})=\big(({a}_1,\dots,{a}_\ell),(n_1,\dots,n_\ell)\big)\in\cH_{\ell,n}(I)$ and set $N:=\mathcal{O}(\pmb{{a}},\pmb{n})$. Given a function $f\in \C^N(I)$, the problem of classical Hermite interpolation is to find a polynomial $P$ of degree at most $n-1$ such that
\Eq{GHI0}{
  f^{(j)}({a}_i)=P^{(j)}({a}_i)
  \qquad(i\in\{1,\dots,\ell\},\,j\in\{0,\dots,n_i-1\})
}
holds. 

Define, for $\alpha\in\{1,\dots,\ell\},\,\beta\in\{0,\dots,n_i-1\}$, the polynomial $H_{\alpha,\beta}$ as follows:
\Eq{Hab}{
  H_{\alpha,\beta}(x)
  :=\sum_{k=0}^{n_\alpha-1-\beta}\frac{\omega(x)}{\beta!k!(x-a_\alpha)^{n_\alpha-\beta-k}}\bigg[\frac{(x-a_\alpha)^{n_\alpha}}{\omega(x)}\bigg]^{(k)}_{x=a_\alpha},
}
where
\Eq{*}{
  \omega(x):=\prod_{i=1}^\ell(x-a_i)^{n_i}.
}
Then, for all $i,\alpha\in\{1,\dots,\ell\}$ and $j,\beta\in\{0,\dots,n_i-1\}$, one can show that
\Eq{*}{
  H_{\alpha,\beta}^{(j)}(a_i)=\delta_{i,\alpha}\cdot\delta_{j,\beta}.
}
Therefore, the function $P$ defined by the formula 
\Eq{*}{
  P(x):=\sum_{\alpha=1}^\ell\sum_{\beta=0}^{n_\alpha-1}
        f^{(\beta)}(a_\alpha)H_{\alpha,\beta}(x)
}
is the unique polynomial of degree less than or equal to $n-1$, which satisfies condition \eq{GHI0}. For the error of the approximation, the following result is known.

\Thm{*}{Let $\ell,n\in\N$ and $(\pmb{{a}},\pmb{n})=\big(({a}_1,\dots,{a}_\ell),(n_1,\dots,n_\ell)\big)\in\cH_{\ell,n}(I)$. 
Define the subintervals $I_0,I_1,\dots,I_\ell$ of $I$ as follows
\Eq{*}{
I_0:=(-\infty,a_1]\cap I,\qquad
I_r:=(a_r,a_{r+1}] \quad(r\in\{1,\dots,\ell-1\}),\qquad
I_{\ell}:=(a_\ell,\infty)\cap I
}
and, for $t\in I$, let $r(t)$ denote the unique index in $\{0,1,\dots,\ell\}$ such that $t\in I_{r(t)}$. Then, for every function $f\in\C^n(I)$ and for all $x\in I$,
\Eq{Cl}{
  f(x)=\sum_{\alpha=1}^\ell\sum_{\beta=0}^{n_\alpha-1}
    f^{(\beta)}(a_\alpha)H_{\alpha,\beta}(x)
    +\int_I f^{(n)}(t) G(x,t) dt,
}
where, for $\alpha\in\{1,\dots,\ell\},\,\beta\in\{0,\dots,n_i-1\}$, the polynomial $H_{\alpha,\beta}$ is defined by \eq{Hab} and the function $G:I\times I\to\R$ is given by 
\Eq{*}{
  G(x,t)
  :=\begin{cases}
  \displaystyle\sum_{i=1}^{r(t)}\sum_{j=0}^{n_j-1} 
  \dfrac{(a_i-t)^{n-j-1}}{(n-j-1)!}H_{i,j}(x) 
  &\mbox{if } t\leq x,\\[5mm]
  \displaystyle-\sum_{i=r(t)+1}^\ell\sum_{j=0}^{n_j-1}\dfrac{(a_i-t)^{n-j-1}}{(n-j-1)!}H_{i,j}(x) 
  &\mbox{if } x<t.
  \end{cases}
}
}
The proof of this theorem can be found, for instance, in the book \cite{AgaWon93}. Several generalizations and applications of this result can be found in the literature, e.g.,  see the papers \cite{AraPecVuk16}, \cite{GodMil81}, \cite{KhaPec20}, \cite{Wan10}.

It seems to be a natural problem to obtain similar interpolations in terms of linear combinations of a given finite set of functions, in particular, in terms of exponential polynomials, which span the kernel of a linear differential operator with constant coefficients. 

The rest of this paper is organized as follows. In Section 2, we recall from \cite{AliPal22,AliPal24} the description of the so-called characteristic element of the kernel of a linear differential operator and establish an addition formula for the elements of the kernel. Several particular cases of this formula are mentioned, e.g., the addition formula for trigonometric and hyperbolic functions and the binomial theorem. In Section 3, we define interpolational systems and describe a generalized Hermite-type interpolation in terms of an $n$-dimensional subspaces of the space of $n$ times continuously differentiable functions. We also introduce the notion of the standard base of such a subspace with respect to the given interpolational system. The standard base is determined in several particular cases. The main result, which is presented in Section 4, establishes an error term in integral form for the Hermite-type interpolation introduced in the previous section in the case when the $n$-dimensional subspace is the kernel of an $n$th order linear differential operator with constant coefficients. Then, several corollaries are deduced illustrating the applicability of our main result.

\section{Auxiliary results on linear differential equations}


For $c=(c_0,\dots,c_n)\in\K^{n+1}$ with $c_n=1$, let $n$th-order linear differential operator $D_c\colon\C^n(I)\to\C(I)$ be defined by the formula
\Eq{D}{
 D_c(f):=c_n f^{(n)}+\dots+c_1 f'+c_0 f \qquad(f\in\C^n(I)).
}
Let $P_c$ denote the \emph{characteristic polynomial of $D_c$}, which is given by
\Eq{P}{
 P_c(\lambda):=c_n \lambda^{n}+\dots+c_1 \lambda+c_0 \qquad(\lambda\in\CC).
}
Let $\lambda_1,\dots,\lambda_k\in\CC$ be pairwise distinct roots $P_c$ with multiplicities $m_1,\dots,m_k\in\N$, respectively. From the theory of linear differential equations, it is well-known that null space of the differential operator $D_c$, denoted as $\ker(D_c)$, is an $n$-dimensional linear subspace of $\C^n(\R)$ consisting of \emph{exponential polynomials}, which is spanned by the following family of \emph{exponential monomials}:
\Eq{*}{
  \R\ni t\mapsto t^j\exp(\lambda_i t),
  \qquad(i\in\{1,\dots,k\},\,j\in\{0,\dots,m_i-1\}).
}
It is also elementary to see that $\ker(D_c)$ is closed with respect to differentiation and translation, i.e., if $\omega\in\ker(D_c)$ and $a\in\R$, we have that $\omega',\tau_a\omega\in\ker(D_c)$, where $(\tau_a\omega)(t):=\omega(t+a)$, $t\in\R$. 

Let $\omega_c\in\C^n(\R)$ denote the unique solution of the initial value problem
\Eq{IV}{
  D_c(\omega_c)=0,\qquad \omega_c^{(\ell)}(0)=\delta_{\ell,n-1} \quad(\ell\in\{0,\dots,n-1\}).
} 
The function $\omega_c$ will be called the \emph{characteristic solution of the differential equation $D_c(\omega)=0$} or the \emph{characteristic element of the kernel of $D_c$}. One can see that if $c\in\R^{n+1}$, then $\omega_c$ is real-valued. The following result, which offers an explicit formula for $\omega_c$, was obtained in \cite{AliPal22}.

\Lem{P}{Let $n\in\N$, $c=(c_0,\dots,c_n)\in\CC^{n+1}$ with $c_n=1$, and let $\omega_c$ be characteristic element of the kernel of $D_c$. Let $\lambda_1,\dots,\lambda_k\in\CC$ be pairwise distinct roots of the characteristic polynomial $P_c$ with multiplicities $m_1,\dots,m_k\in\N$, respectively. Define the polynomials $P_1,\dots,P_k:\CC\to\CC$ by
\Eq{*}{
    P_i(\lambda):=\prod_{\ell\in\{1,\dots,k\}\setminus\{i\}}(\lambda-\lambda_\ell)^{m_\ell} \qquad(i\in\{1,\dots,k\}).
}
Then, for all $t\in\R$,
\Eq{*}{
  \omega_c(t)
  =\sum_{i=1}^k\sum_{j=0}^{m_i-1}\frac{(P_i^{-1})^{(m_i-1-j)}(\lambda_i)}{(m_i-1-j)!j!}\cdot t^j\exp(\lambda_i t).
}}

We shall also need the following result from the paper \cite{AliPal22} which was not stated explicitly therein, however is formulated and proved in the proof of Theorem 3.1.

\Lem{PP}{Let $n\in\N$, $c=(c_0,\dots,c_n)\in\CC^{n+1}$ with $c_n=1$, and let $\omega_c$ denote the characteristic solution of the differential operator $D_c$. Then, for all $j,\beta\in\{0,\dots,n-1\}$,
\Eq{*}{
  \delta_{j,\beta}=\sum_{i=0}^{n-j-1}
  c_{i+j+1}\cdot\omega_c^{(\beta+i)}(0).
}}

Using this lemma, we can obtain the following addition theorem for the elements of the kernel of the differential operator $D_c$.

\Thm{AF}{Let $n\in\N$, $c=(c_0,\dots,c_n)\in\CC^{n+1}$ with $c_n=1$, and let $\omega_c$ denote the characteristic solution of the differential operator $D_c$. Then, for all $\omega\in\ker(D_c)$ and for all $u,v\in\R$, 
\Eq{uv}{
  \omega(u+v)=\sum_{j=0}^{n-1}\sum_{i=0}^{n-1-j}c_{i+j+1}\omega_c^{(i)}(u)\omega^{(j)}(v).
}}

In view of the formula \eq{uv}, we can see that each member of the kernel of $D_c$ satisfies a Levi--Civita-type functional equation. 

\begin{proof} Let $v\in\R$ be fixed. As functions of the variable $u$, both sides of the formula \eq{uv} are elements of the kernel of $D_c$, i.e., the functions $\varphi,\psi:\R\to\CC$ defined by
\Eq{*}{
  \varphi(u):=\omega(u+v),\qquad
  \psi(u):=\sum_{j=0}^{n-1}\sum_{i=0}^{n-1-j}c_{i+j+1}\omega_c^{(i)}(u)\omega^{(j)}(v)
}
belong to $\ker(D_c)$. On the other hand, for all $\beta\in\{0,\dots,n-1\}$,
\Eq{*}{
  \varphi^{(\beta)}(0):=\omega^{(\beta)}(v)
}
and, according to \lem{PP}, we can obtain that
\Eq{*}{
  \psi^{(\beta)}(0)=\sum_{j=0}^{n-1}\bigg(\sum_{i=0}^{n-1-j}c_{i+j+1}\omega_c^{(\beta+i)}(0)\bigg)\omega^{(j)}(v)
  =\sum_{j=0}^{n-1}\delta_{j,\beta}\omega^{(j)}(v)
  =\omega^{(\beta)}(v).
}
Therefore, for all $\beta\in\{0,\dots,n-1\}$, we have that 
$\varphi^{(\beta)}(0)=\omega^{(\beta)}(v)=\psi^{(\beta)}(0)$. By the uniqueness of the solutions of the linear differential equation $D_c(f)=0$ with a given initial value condition at $u=0$, it follows that $\varphi$ and $\psi$ are identical functions. This proves that the equality \eq{uv} holds for all $u,v\in\R$. 
\end{proof}

\Exa{AFH}{Let $n=2$ and $c:=(c_0,c_1,c_2):=(-1,0,1)$. Then the differential operator $D_c$ is given as $D_c(f)=f''-f$. The kernel of $D_c$ is spanned by the functions $\cosh$ and $\sinh$ and the characteristic solution $\omega_c:\R\to\R$ is given by
\Eq{1Ho}{
  \omega_c=\sinh.
}
According to \thm{AF}, for all $u,v\in\R$, we have that
\Eq{*}{
  \cosh(u+v)&=\sum_{j=0}^1\sum_{i=0}^{1-j}c_{i+j+1}\sinh^{(i)}(u)\cosh^{(j)}(v)\\
  &=c_{1}\sinh^{(0)}(u)\cosh^{(0)}(v)+c_{2}\sinh^{(1)}(u)\cosh^{(0)}(v)+c_{2}\sinh^{(0)}(u)\cosh^{(1)}(v).
  \\
  &=\cosh(u)\cosh(v)+\sinh(u)\sinh(v),
}
and
\Eq{*}{
  \sinh(u+v)&=\sum_{j=0}^1\sum_{i=0}^{1-j}c_{i+j+1}\sinh^{(i)}(u)\sinh^{(j)}(v)\\
  &=c_{1}\sinh^{(0)}(u)\sinh^{(0)}(v)+c_{2}\sinh^{(1)}(u)\sinh^{(0)}(v)+c_{2}\sinh^{(0)}(u)\sinh^{(1)}(v).
  \\
  &=\cosh(u)\sinh(v)+\sinh(u)\cosh(v).
}
These equalities are the standard addition formulae for hyperbolic functions.
}

\Exa{AFT}{Let $n=2$ and $c:=(c_0,c_1,c_2):=(1,0,1)$. Then the differential operator $D_c$ is given as $D_c(f)=f''+f$. The kernel of $D_c$ is spanned by the functions $\cos$ and $\sin$ and the characteristic solution $\omega_c:\R\to\R$ is given by
\Eq{1To}{
  \omega_c=\sin.
}
According to \thm{AF}, for all $u,v\in\R$, we have that
\Eq{*}{
  \cos(u+v)&=\sum_{j=0}^1\sum_{i=0}^{1-j}c_{i+j+1}\sin^{(i)}(u)\cos^{(j)}(v)\\
  &=c_{1}\sin^{(0)}(u)\cos^{(0)}(v)+c_{2}\sin^{(1)}(u)\cos^{(0)}(v)+c_{2}\sin^{(0)}(u)\cos^{(1)}(v).
  \\
  &=\cos(u)\cos(v)-\sin(u)\sin(v),
}
and
\Eq{*}{
  \sin(u+v)&=\sum_{j=0}^1\sum_{i=0}^{1-j}c_{i+j+1}\sin^{(i)}(u)\sin^{(j)}(v)\\
  &=c_{1}\sin^{(0)}(u)\sin^{(0)}(v)+c_{2}\sin^{(1)}(u)\sin^{(0)}(v)+c_{2}\sin^{(0)}(u)\sin^{(1)}(v).
  \\
  &=\cos(u)\sin(v)+\sin(u)\cos(v).
}
These equalities are the standard addition formulae for trigonometric functions.
}

\Exa{AFP}{Let $n\in\N$ and $c:=(c_0,\dots,c_{n-1},c_n):=(0,\dots,0,1)$. Then the differential operator $D_c$ is given as $D_c(f)=f^{(n)}$. The kernel of $D_c$ is spanned by the polynomials of degree less than or equal to $n-1$ and the characteristic solution $\omega_c:\R\to\R$ is given by
\Eq{3CE}{
  \omega_c(t)=\frac{t^{n-1}}{(n-1)!} \qquad (t\in\R).
}
According to \thm{AF}, for all $k\in\{0,\dots,n-1\}$ and $u,v\in\R$, we obtain that
\Eq{*}{
  (u+v)^k&=\sum_{j=0}^{n-1}\sum_{i=0}^{\min(k,n-1-j)}c_{i+j+1}\frac{k!u^{k-i}}{(k-i)!}\frac{v^{n-j-1}}{(n-j-1)!}\\
  &=\sum_{j=n-1-k}^{n-1}\frac{k!u^{k-(n-1-j)}}{(k-(n-1-j))!}\frac{v^{n-j-1}}{(n-j-1)!}\\
  &=\sum_{i=0}^k\frac{k!u^{k-i}}{(k-i)!}\frac{v^i}{i!}
  =\sum_{i=0}^k\binom{k}{i}u^{k-i}v^i,
}
which shows that the binomial theorem is also a particular case of \thm{AF}.
}

\section{Hermite-type interpolation with exponential polynomials}

Assume that $\Omega$ is an $n$-dimensional linear subspace of $\C^n(I)$. We are going to construct Hermite-type interpolation for the elements of $\C^n(I)$ in terms of the elements of $\Omega$. To describe such interpolations, for $\ell,n\in\N$, consider the set
\Eq{*}{
  \cH_{\ell,n}(I):=\Big\{(\pmb{{a}},\pmb{n})\mid\, 
  & \pmb{{a}}=({a}_1,\dots,{a}_\ell)\in I^\ell \mbox{ with } {a}_1<\dots<{a}_\ell  \mbox{ and}\\
  & \pmb{n}=(n_1,\dots,n_\ell)\in\N^\ell \mbox{ with } n=n_1+\dots+n_\ell\Big\}.
}
Given an interpolational system $(\pmb{{a}},\pmb{n})=\big(({a}_1,\dots,{a}_\ell),(n_1,\dots,n_\ell)\big)\in\cH_{\ell,n}(I)$, the number $n$ will be called its \emph{dimension}. The order $\mathcal{O}(\pmb{{a}},\pmb{n})$ of $(\pmb{{a}},\pmb{n})$ is defined to be number $\max(n_1,\dots,n_\ell)-1$.

Let $(\pmb{{a}},\pmb{n})=\big(({a}_1,\dots,{a}_\ell),(n_1,\dots,n_\ell)\big)\in\cH_{\ell,n}(I)$ and set $N:=\mathcal{O}(\pmb{{a}},\pmb{n})$. Given an $n$-dimensional subspace $\Omega$ of $\C^N(I)$ and function $f\in \C^N(I)$, the problem is to find an element $\omega\in\Omega$ such that
\Eq{GHI}{
  f^{(j)}({a}_i)=\omega^{(j)}({a}_i)
  \qquad(i\in\{1,\dots,\ell\},\,j\in\{0,\dots,n_i-1\})
}
hold. The following lemma offers a condition for the unique solvability of the above interpolation problem.

\Lem{1}{Let $\ell,n\in\N$, $(\pmb{{a}},\pmb{n})=\big(({a}_1,\dots,{a}_\ell),(n_1,\dots,n_\ell)\big)\in\cH_{\ell,n}(I)$ and set $N:=\mathcal{O}(\pmb{{a}},\pmb{n})$. Let $\Omega$ be an $n$-dimensional linear subspace of $\C^N(I)$ and let $\omega_1,\dots,\omega_n\in\Omega$ be a linear base for $\Omega$. Assume that
\Eq{NDC}{
  W_{\pmb{{a}}}^{\pmb{n}}
  (\omega_1,\dots,\omega_n):=
  \begin{vmatrix}
  \omega_1({a}_1) & \dots & \omega_1^{(n_1-1)}({a}_1) & \quad \dots \quad & \omega_1({a}_\ell) & \dots & \omega_1^{(n_\ell-1)}({a}_\ell) \\
  \vdots & \ddots & \vdots & \ddots & \vdots & \ddots &\vdots \\
  \omega_n({a}_1) & \dots & \omega_n^{(n_1-1)}({a}_1) & \dots & \omega_n({a}_\ell) & \dots & \omega_n^{(n_\ell-1)}({a}_\ell) \\  
  \end{vmatrix}\neq0.
}
holds. Then, for any $f\in \C^N(I)$, the system of linear equations 
\Eq{SE}{
  f^{(j)}({a}_i)=\sum_{\alpha=1}^n\mu_\alpha\omega_\alpha^{(j)}({a}_i)
  \qquad(i\in\{1,\dots,\ell\},\,j\in\{0,\dots,n_i-1\}).
}
has a unique solution $\mu_1,\dots,\mu_n\in\CC$ and hence the function $\omega=\sum_{\alpha=1}^n\mu_\alpha\omega_\alpha$ satisfies the interpolation condition \eq{GHI}. Furthermore, if $f\in\Omega$, then $f=\omega$.}

\begin{proof} Let $f\in \C^N(I)$. Then, the base determinant of the system of linear equations 
\eq{SE} is non-vanishing according to condition \eq{NDC}. This ensures the unique solvability of \eq{SE} with respect to $\mu_1,\dots,\mu_n\in\CC$ follows. The equalities in \eq{SE} also show that $\omega=\sum_{\alpha=1}^n\mu_\alpha\omega_\alpha$ satisfies the interpolation condition \eq{GHI}.

If $f$ belongs to $\Omega$, then $f$ is a linear combination of the functions $\omega_1,\dots,\omega_n$. Therefore, there exist $\nu_1,\dots,\nu_n\in\CC$ such that $f=\sum_{\alpha=1}^n\nu_\alpha \omega_\alpha$. Then, obviously,
\Eq{*}{
  f^{(j)}({a}_i)=\sum_{\alpha=1}^n\nu_\alpha\omega_\alpha^{(j)}({a}_i)
  \qquad(i\in\{1,\dots,\ell\},\,j\in\{0,\dots,n_i-1\}).
}
Combining this system of equations with \eq{SE}, we can see that
\Eq{*}{
  0=\sum_{\alpha=1}^n(\nu_\alpha-\mu_\alpha)\omega_\alpha^{(j)}({a}_i)
  \qquad(i\in\{1,\dots,\ell\},\,j\in\{0,\dots,n_i-1\}).
}
Due to condition \eq{NDC}, it follows that $(\nu_1,\dots,\nu_n)=(\mu_1,\dots,\mu_n)$, and hence $f=\sum_{\alpha=1}^n\mu_\alpha \omega_\alpha$ must be valid.
\end{proof}

\Lem{2}{Let $\ell,n\in\N$, $(\pmb{{a}},\pmb{n})=\big(({a}_1,\dots,{a}_\ell),(n_1,\dots,n_\ell)\big)\in\cH_{\ell,n}(I)$ and set $N:=\mathcal{O}(\pmb{{a}},\pmb{n})$. Let $\Omega$ be an $n$-dimensional linear subspace of $\C^N(I)$ and let $\omega_1,\dots,\omega_n\in\Omega$ be a linear base for $\Omega$. Assume that \eq{NDC} holds.
Then, for any $\alpha\in\{1,\dots,\ell\}$ and $\beta\in\{0,\dots,n_\alpha-1\}$, there exists a unique element $\chi_{\alpha,\beta}\in\Omega$ which satisfies the following conditions
\Eq{chi}{
  \chi_{\alpha,\beta}^{(j)}({a}_i)
  =\delta_{i,\alpha}\delta_{j,\beta}
  \qquad(i\in\{1,\dots,\ell\},\,j\in\{0,\dots,n_i-1\}).
}
Furthermore, the system $\{\chi_{\alpha,\beta}\colon \alpha\in\{1,\dots,\ell\},\beta\in\{0,\dots,n_\alpha-1\}\}$ forms a base of $\Omega$ and, for any $f\in\C^N(I)$, 
\Eq{om}{
  \omega
  =\sum_{\alpha=1}^\ell\sum_{\beta=0}^{n_\alpha-1} f^{(\beta)}({a}_\alpha)\chi_{\alpha,\beta}
}
is the unique element of $\Omega$ which satisfies the interpolation condition \eq{GHI}. In addition, if $f\in\Omega$, then $f=\omega$.}

\begin{proof} Let $\alpha\in\{1,\dots,\ell\}$ and $\beta\in\{0,\dots,n_\alpha-1\}$ be fixed. To prove the existence of $\chi_{\alpha,\beta}\in\Omega$ satisfying \eq{chi}, we search for  
$\chi_{\alpha,\beta}$ in the form $\chi_{\alpha,\beta}=\sum_{k=1}^n\mu_k\omega_k$ with some constants $\mu_1,\dots \mu_n\in\CC$. Then\eq{chi} reduces to the following linear system of equations $\mu_1,\dots \mu_n\in\CC$:
\Eq{*}{
  \delta_{i,\alpha}\delta_{j,\beta}
  =\sum_{k=1}^n\mu_k\omega_k^{(j)}({a}_i)
  \qquad(i\in\{1,\dots,\ell\},\,j\in\{0,\dots,n_i-1\}).
}
The base determinant of this system is nonvanishing according to our assumption \eq{NDC}, therefore this system of equations is uniquely solvable for $\mu_1,\dots \mu_n\in\CC$ and this proves the existence of $\chi_{\alpha,\beta}\in\Omega$ satisfying \eq{chi}.

To point out that the system $\{\chi_{\alpha,\beta}\colon \alpha\in\{1,\dots,\ell\},\beta\in\{0,\dots,n_\alpha-1\}\}$ forms a base of $\Omega$, it suffices to show its linear independence. Assume that, for some constants $c_{\alpha,\beta}\in\CC$, $\alpha\in\{1,\dots,\ell\}$, $\beta\in\{0,\dots,n_\alpha-1\}$, we have
\Eq{lid}{
  0=\sum_{\alpha=1}^\ell\sum_{\beta=0}^{n_\alpha-1} c_{\alpha,\beta}\chi_{\alpha,\beta}.
}
Then, for all $i\in\{1,\dots,\ell\}$ and $j\in\{0,\dots,n_i-1\}$,
\Eq{*}{
  0=\bigg(\sum_{\alpha=1}^\ell\sum_{\beta=0}^{n_\alpha-1} c_{\alpha,\beta}
  \chi_{\alpha,\beta}\bigg)^{(j)}({a}_i)
  =\sum_{\alpha=1}^\ell\sum_{\beta=0}^{n_\alpha-1} c_{\alpha,\beta}
  \chi_{\alpha,\beta}^{(j)}({a}_i)
  =\sum_{\alpha=1}^\ell\sum_{\beta=0}^{n_\alpha-1}c_{\alpha,\beta}
  \delta_{i,\alpha}\delta_{j,\beta}
  =c_{i,j}.
}
Hence \eq{lid} can only hold when $c_{\alpha,\beta}=0$ for all $\alpha\in\{1,\dots,\ell\}$, $\beta\in\{0,\dots,n_\alpha-1\}$. This shows that the system $\{\chi_{\alpha,\beta}\colon \alpha\in\{1,\dots,\ell\},\beta\in\{0,\dots,n_\alpha-1\}\}$ is linearly independent, indeed.

Being a linear combination of elements of $\Omega$, it is obvious that the function \eq{om} belongs to $\Omega$. Furthermore, for all $i\in\{1,\dots,\ell\}$ and $j\in\{0,\dots,n_i-1\}$, we also have that
\Eq{*}{
  \omega^{(j)}({a}_i)
  =\bigg(\sum_{\alpha=1}^\ell\sum_{\beta=0}^{n_\alpha-1} f^{(\beta)}({a}_\alpha)\chi_{\alpha,\beta}\bigg)^{(j)}({a}_i)
  &=\sum_{\alpha=1}^\ell\sum_{\beta=0}^{n_\alpha-1} f^{(\beta)}({a}_\alpha)\chi_{\alpha,\beta}^{(j)}({a}_i)\\
  &=\sum_{\alpha=1}^\ell\sum_{\beta=0}^{n_\alpha-1}f^{(\beta)}({a}_\alpha)\delta_{i,\alpha}\delta_{j,\beta}
  =f^{(j)}({a}_i).
}
Therefore, $\omega$ is the unique element of $\Omega$ which satisfies the interpolation condition \eq{GHI}. 

If, in  addition $f\in\Omega$, then, according to the last statement of \lem{1}, for its interpolating function $\omega\in\Omega$, we have that $f=\omega$ holds. This implies the last statement of this lemma.  
\end{proof}

The system $\{\chi_{\alpha,\beta}\colon \alpha\in\{1,\dots,\ell\},\beta\in\{0,\dots,n_\alpha-1\}\}$ constructed in the above lemma will be called the \emph{standard base of $\Omega$ related to the interpolational system $(\pmb{{a}},\pmb{n})\in\cH_{\ell,n}(I)$}. 

\medskip

In what follows, we consider the case when $(\omega_1,\dots,\omega_n):\R\to\CC^n$ forms a fundamental system of solutions of the $n$th-order homogeneous linear differential equation $D_c(f)=0$, where $n\in\N$, $c=(c_0,\dots,c_n)\in\CC^{n+1}$ with $c_n=1$, and the differential operator $D_c:\C^n(I)\to\C(I)$ is defined by \eq{D}. Let $\lambda_1,\dots,\lambda_k\in\CC$ denote the pairwise distinct roots of the characteristic polynomial $P_c$ with multiplicities $m_1,\dots,m_k\in\N$, respectively. Then, according to the theory of higher-order linear differential equation, a fundamental system of solutions $(\omega_1,\dots,\omega_n):\R\to\CC^n$ of the differential equation $D_c(f)=0$ can be given by
\Eq{*}{
  (\omega_1,\dots,\omega_n)(t)
  :=\big(\exp(\lambda_1 t),\dots,t^{m_1-1}\exp(\lambda_1 t),\dots,\exp(\lambda_k t),\dots,t^{m_k-1}\exp(\lambda_k t)\big).
}
In the first two examples, we describe the settings related to standard Lagrange and Taylor interpolations, respectively.

\Exa{3}{Let $n\in\N$, let $\big(({a}_1,\dots,{a}_n),(1,\dots,1)\big)\in\cH_{n,n}(I)$ and  $c:=(0,\dots,0,1)\in\R^{n+1}$. Then the differential operator $D_c\colon\C^n(I)\to\C(I)$ is given by $D_c(f):=f^{(n)}$. For $i\in\{1,\dots,n\}$, define 
\Eq{3SB}{
  \chi_{i,0}(t)
  =\prod_{j\in\{1,\dots,n\}\setminus\{i\}}
    \frac{t-{a}_j}{{a}_i-{a}_j} \qquad(t\in\R).
}
These polynomials of degree $n-1$ and hence they belong to the null space of $D_c$. (They are usually called the basic Lagrange interpolational polynomials.) It is easy to see that
\Eq{*}{
  \chi_{i,0}({a}_\alpha)=\delta_{i,\alpha}\qquad (i,\alpha\in\{1,\dots,n\}).
}
Consequently, $\{\chi_{1,0},\dots,\chi_{n,0}\}$ is the standard base of the null space of $D_c$ related to the interpolational system $\big(({a}_1,\dots,{a}_n),(1,\dots,1)\big)$. The characteristic element $\omega_c:\R\to\R$ of $\ker(D_c)$ now has the form \eq{3CE}.}

\Exa{4}{Let $n\in\N$, let $({a},n)\in\cH_{n,1}(I)$ and let $c:=(0,\dots,0,1)\in\R^{n+1}$. Then the differential operator $D_c\colon\C^n(I)\to\C(I)$ be given by $D(f):=f^{(n)}$. For $j\in\{0,\dots,n-1\}$, define 
\Eq{4SB}{
  \chi_{1,j}(t)
  =\frac{(t-{a})^j}{j!} \qquad(t\in\R).
}
These polynomials of at most degree $n-1$ and hence they belong to the null space of $D$. (In fact, they are the basic Taylor interpolational polynomials at the point ${a}$.) It is easy to see that
\Eq{*}{
  \chi_{1,j}^{(\beta)}({a})=\delta_{j,\beta}\qquad (j,\beta\in\{0,\dots,n-1\}).
}
Therefore, $\{\chi_{1,0},\dots,\chi_{1,n-1}\}$ is the standard base of the kernel of $D_c$ related to the interpolational system $({a},n)$. The characteristic element $\omega_c:\R\to\R$ of $\ker(D_c)$ is given by \eq{3CE}.}

\Exa{4+}{Let $n\in\N$, let $({a},n)\in\cH_{n,1}(I)$ and let $c:=(c_0,\dots,c_n)\in\CC^{n+1}$ and define the differential operator $D_c\colon\C^n(I)\to\C(I)$ be given by \eq{D}. The characteristic element $\omega_c:\R\to\R$ of $\ker(D_c)$ can be computed with the help of \lem{P}. For $j\in\{0,\dots,n-1\}$, define 
\Eq{4SB+}{
  \chi_{1,j}(t)
  :=\sum_{i=0}^{n-j-1} c_{i+j+1}\cdot\omega_c^{(i)}(t-{a})
  \qquad (t\in\R.
}
The kernel of $D_c$ is closed with respect to differentiation and translation of the argument, therefore, the each of the functions $\chi_{1,0},\dots,\chi_{1,n-1}$ belongs to $\ker(D_c)$. On the other hand, according to \lem{PP}, we have that
\Eq{*}{
  \chi_{1,j}^{(\beta)}({a})=\delta_{j,\beta}\qquad (j,\beta\in\{0,\dots,n-1\}).
}
Consequently, $\{\chi_{1,0},\dots,\chi_{1,n-1}\}$ is the standard base of the kernel of $D_c$ related to the interpolational system $({a},n)$. }

\Exa{1H}{Let $\big(({a}_1,{a}_2),(1,1)\big)\in\cH_{2,2}(I)$ and let $c:=(-1,0,1)$. Then ${a}_1\neq{a}_2$ and the second-order differential operator $D_c\colon\C^2(I)\to\C(I)$ is given as $D_c(f):=f''-f$. Define
\Eq{1H}{
  \chi_{1,0}(t):=\frac{\sinh(t-{a}_2)}{\sinh({a}_1-{a}_2)}
  \qquad\mbox{and}\qquad
  \chi_{2,0}(t):=\frac{\sinh(t-{a}_1)}{\sinh({a}_2-{a}_1)} \qquad(t\in\R).
}
These functions are linear combinations of the hyperbolic functions $\sinh$ and $\cosh$, therefore they belong to the null space of $D_c$. On the other hand, 
\Eq{SBC}{
  \chi_{i,0}({a}_\alpha)=\delta_{i,\alpha} \qquad (i,\alpha\in\{1,2\}).
}
Therefore, $\{\chi_{1,0},\chi_{2,0}\}$ is the standard base of the null space of $D_c$ related to $\big(({a}_1,{a}_2),(1,1)\big)$. One can also see that the characteristic solution $\omega_c:\R\to\R$ is given by \eq{1Ho}.}

\Exa{1T}{Let $\big(({a}_1,{a}_2),(1,1)\big)\in\cH_{2,2}(I)$ and assume that $\sin({a}_1-{a}_2)\neq0$. Let $c:=(1,0,1)$. Then the second-order differential operator $D_c\colon\C^2(I)\to\C(I)$ is given as $D_c(f):=f''+f$. Define
\Eq{1T}{
  \chi_{1,0}(t):=\frac{\sin(t-{a}_2)}{\sin({a}_1-{a}_2)}
  \qquad\mbox{and}\qquad
  \chi_{2,0}(t):=\frac{\sin(t-{a}_1)}{\sin({a}_2-{a}_1)} \qquad(t\in\R).
}
These functions are linear combinations of the trigonometric functions $\sin$ and $\cos$, thus they belong to the null space of $D$. One can also see that \eq{SBC} holds, thus, $\{\chi_{1,0},\chi_{2,0}\}$ forms the standard base of the kernel of $D$ related to $\big(({a}_1,{a}_2),(1,1)\big)$. The characteristic element $\omega_c:\R\to\R$ of the null space of $D_c$ is now given by \eq{1To}.}

\Exa{2HT}{Let $\big(({a}_1,{a}_2),(2,2)\big)\in\cH_{2,4}(I)$ and assume that $\cosh({a}_1-{a}_2)\cos({a}_1-{a}_2)\neq 1$. Let $c:=(-1,0,0,0,1)$. Then the $4$th-order differential operator $D_c\colon\C^4(I)\to\C(I)$ is defined by $D_c(f):=f^{(4)}-f$ and the functions 
\Eq{*}{
\omega_1(t):=\sinh(t-{a}_2),\qquad 
\omega_2(t):=\cosh(t-{a}_2),\qquad 
\omega_3(t):=\sin(t-{a}_2),\qquad 
\omega_4(t):=\cos(t-{a}_2)
}
form a base for the null space of $D_c$. Then
\Eq{*}{
 W_{{a}_1,{a}_2}^{2,2}(\omega_1,\dots,\omega_4)
 &=\begin{vmatrix}
  \sinh({a}_1-{a}_2) & \cosh({a}_1-{a}_2) & 0 & 1 \\
  \cosh({a}_1-{a}_2) & \sinh({a}_1-{a}_2) & 1 & 0 \\
  \sin({a}_1-{a}_2) & \cos({a}_1-{a}_2) & 0 & 1 \\
  \cos({a}_1-{a}_2) & -\sin({a}_1-{a}_2) & 1 & 0 \\
  \end{vmatrix}
  =2-2\cosh({a}_1-{a}_2)\cos({a}_1-{a}_2),
}
which is nonzero according to our assumption. For $t\in\R$, define 
\comment{
\Eq{*}{
  \chi_{1,0}(t)&:=\frac{1}{W_{({a}_1,{a}_2)}^{(2,2)}(\omega_1,\dots,\omega_4)}
  \begin{vmatrix}
  \sinh(t-{a}_2) & \cosh({a}_1-{a}_2) & 0 & 1 \\
  \cosh(t-{a}_2) & \sinh({a}_1-{a}_2) & 1 & 0 \\
  \sin(t-{a}_2) & \cos({a}_1-{a}_2) & 0 & 1 \\
  \cos(t-{a}_2) & -\sin({a}_1-{a}_2) & 1 & 0 \\
  \end{vmatrix},\\
  \chi_{1,1}(t)&:=\frac{1}{W_{({a}_1,{a}_2)}^{(2,2)}(\omega_1,\dots,\omega_4)}
  \begin{vmatrix}
  \sinh({a}_1-{a}_2) & \sinh(t-{a}_2) & 0 & 1 \\
  \cosh({a}_1-{a}_2) & \cosh(t-{a}_2) & 1 & 0 \\
  \sin({a}_1-{a}_2) & \sin(t-{a}_2) & 0 & 1 \\
  \cos({a}_1-{a}_2) & \cos(t-{a}_2) & 1 & 0 \\
  \end{vmatrix},\\
  \chi_{2,0}(t)&:=\frac{1}{W_{({a}_1,{a}_2)}^{(2,2)}(\omega_1,\dots,\omega_4)}
  \begin{vmatrix}
  \sinh({a}_1-{a}_2) & \cosh({a}_1-{a}_2) & \sinh(t-{a}_2) & 1 \\
  \cosh({a}_1-{a}_2) & \sinh({a}_1-{a}_2) & \cosh(t-{a}_2) & 0 \\
  \sin({a}_1-{a}_2) & \cos({a}_1-{a}_2) & \sin(t-{a}_2) & 1 \\
  \cos({a}_1-{a}_2) & -\sin({a}_1-{a}_2) & \cos(t-{a}_2) & 0 \\
  \end{vmatrix},\\
  \chi_{2,1}(t)&:=\frac{1}{W_{({a}_1,{a}_2)}^{(2,2)}(\omega_1,\dots,\omega_4)}
  \begin{vmatrix}
  \sinh({a}_1-{a}_2) & \cosh({a}_1-{a}_2) & 0 & \sinh(t-{a}_2) \\
  \cosh({a}_1-{a}_2) & \sinh({a}_1-{a}_2) & 1 & \cosh(t-{a}_2) \\
  \sin({a}_1-{a}_2) & \cos({a}_1-{a}_2) & 0 & \sin(t-{a}_2) \\
  \cos({a}_1-{a}_2) & -\sin({a}_1-{a}_2) & 1 & \cos(t-{a}_2) \\
  \end{vmatrix}.
}
Then, after a short computation, one gets}
\begin{equation}\begin{aligned}\label{E2HT}
  \chi_{1,0}(t)&:=\frac{1}{W_{{a}_1,{a}_2}^{2,2}(\omega_1,\dots,\omega_4)}\begin{vmatrix}
   \cosh({a}_1-{a}_2)-\cos({a}_1-{a}_2) 
   & \sinh(t-{a}_2)-\sin(t-{a}_2) \\
   \sinh({a}_1-{a}_2)+\sin({a}_1-{a}_2) 
   & \cosh(t-{a}_2)-\cos(t-{a}_2)
   \end{vmatrix},\\
  \chi_{1,1}(t)&:=\frac{1}{W_{{a}_1,{a}_2}^{2,2}(\omega_1,\dots,\omega_4)}\begin{vmatrix}
   \cosh({a}_1-{a}_2)-\cos({a}_1-{a}_2) 
   & \cosh(t-{a}_2)-\cos(t-{a}_2) \\
   \sinh({a}_1-{a}_2)-\sin({a}_1-{a}_2) 
   & \sinh(t-{a}_2)-\sin(t-{a}_2)
   \end{vmatrix},\\
  \chi_{2,0}(t)&:=\frac{1}{W_{{a}_1,{a}_2}^{2,2}(\omega_1,\dots,\omega_4)}\begin{vmatrix}
   \cosh({a}_1-{a}_2)-\cos({a}_1-{a}_2) 
   & -\sinh(t-{a}_1)+\sin(t-{a}_1) \\
   \sinh({a}_1-{a}_2)+\sin({a}_1-{a}_2) 
   & \cosh(t-{a}_1)-\cos(t-{a}_1)
   \end{vmatrix},\\
  \chi_{2,1}(t)&:=\frac{1}{W_{{a}_1,{a}_2}^{2,2}(\omega_1,\dots,\omega_4)}\begin{vmatrix}
   \cosh({a}_1-{a}_2)-\cos({a}_1-{a}_2) 
   & -\cosh(t-{a}_1)+\cos(t-{a}_1) \\
   \sinh({a}_1-{a}_2)-\sin({a}_1-{a}_2) 
   & \sinh(t-{a}_1)-\sin(t-{a}_1)
   \end{vmatrix}.
\end{aligned}\end{equation}
Then it is simple to check that
\Eq{*}{
  \chi_{i,j}^{(\beta)}({a}_\alpha)=\delta_{i,\alpha}\delta_{j,\beta}\qquad (i,\alpha\in\{1,2\},\,j,\beta\in\{0,1\})
}
is satisfied. Therefore, $\{\chi_{1,0},\chi_{1,1},\chi_{2,0},\chi_{2,1}\}$ is the standard base of the kernel of $D_c$ related to the interpolational system $\big(({a}_1,{a}_2),(2,2)\big)$. The characteristic element $\omega_c:\R\to\R$ of $\ker(D_c)$ now has the following form
\Eq{2HTo}{
  \omega_c(t)=\frac12\big(\sinh(t)-\sin(t)\big) \qquad(t\in\R).
}}

\Exa{5H}{Let $\big(({a}_1,{a}_2,{a}_3),(1,1,1)\big)\in\cH_{3,3}(I)$ and let $c:=(0,-1,0,1)$. Then the $3$rd-order differential operator $D_c\colon\C^3(I)\to\C(I)$ is given by $D(f):=f'''-f'$, the null space of $D_c$ is spanned by the functions $1,\cosh(t-{a}_3),\sinh(t-{a}_3)$, and we have that
\Eq{*}{
  W_{{a}_1,{a}_2,{a}_3}^{1,1,1}(1,\cosh,\sinh)&=
  \begin{vmatrix}
   1 & 1 & 1 \\
   \cosh({a}_1-{a}_3) & \cosh({a}_2-{a}_3) & 1 \\
   \sinh({a}_1-{a}_3) & \sinh({a}_2-{a}_3) & 0 
  \end{vmatrix}\\
  &=\sinh({a}_3-{a}_2)+\sinh({a}_1-{a}_3)+\sinh({a}_2-{a}_1)\\
  &=4\sinh\Big(\frac{{a}_1-{a}_3}{2}\Big)\sinh\Big(\frac{{a}_2-{a}_1}{2}\Big)\sinh\Big(\frac{{a}_3-{a}_2}{2}\Big)\neq0.
}
Define
\Eq{5H}{
  \chi_{1,0}(t)&:=
  \comment{\frac{1}{W_{{a}_1,{a}_2,{a}_3}^{1,1,1}(1,\cosh,\sinh)}
  \begin{vmatrix}
   1 & 1 & 1\\
   \cosh(t-{a}_3) & \cosh({a}_2-{a}_3) & 1 \\
   \sinh(t-{a}_3) & \sinh({a}_2-{a}_3) & 0
  \end{vmatrix}
  \\
  &=\frac{\sinh\Big(\frac{t-{a}_3}{2}\Big)\sinh\Big(\frac{{a}_2-t}{2}\Big)}{\sinh\Big(\frac{{a}_1-{a}_3}{2}\Big)\sinh\Big(\frac{{a}_2-{a}_1}{2}\Big)}
  =}\frac{\cosh\big(t-\frac{{a}_2+{a}_3}{2}\big)-\cosh\big(\frac{{a}_2-{a}_3}{2}\big)}{\cosh\big({a}_1-\frac{{a}_2+{a}_3}{2}\big)-\cosh\big(\frac{{a}_2-{a}_3}{2}\big)},\\
  \chi_{2,0}(t)&:=
  \frac{\cosh\big(t-\frac{{a}_1+{a}_3}{2}\big)-\cosh\big(\frac{{a}_1-{a}_3}{2}\big)}{\cosh\big({a}_2-\frac{{a}_1+{a}_3}{2}\big)-\cosh\big(\frac{{a}_1-{a}_3}{2}\big)},\\
  \chi_{3,0}(t)&:=
  \frac{\cosh\big(t-\frac{{a}_1+{a}_2}{2}\big)-\cosh\big(\frac{{a}_1-{a}_2}{2}\big)}{\cosh\big({a}_3-\frac{{a}_1+{a}_2}{2}\big)-\cosh\big(\frac{{a}_1-{a}_2}{2}\big)}.
}
One can easily see that $\chi_{1,0},\chi_{2,0},\chi_{3,0}$ form the standard base of the kernel of $D_c$ related to the interpolational system $\big(({a}_1,{a}_2,{a}_3),(1,1,1)\big)$.
The characteristic element $\omega_c:\R\to\R$ of the null space of $D_c$ is now given as
\Eq{5Ho}{
\omega_c(t)=\cosh(t)-1 \qquad(t\in\R).
}
}

\Exa{5T}{Let $\big(({a}_1,{a}_2,{a}_3),(1,1,1)\big)\in\cH_{3,3}(I)$ with $
\sin\big(\frac{{a}_1-{a}_3}{2}\big)\sin\big(\frac{{a}_2-{a}_1}{2}\big)\sin\big(\frac{{a}_3-{a}_2}{2}\big)\neq0$ and let $c:=(0,1,0,1)$. Then the $3$rd-order differential operator $D_c\colon\C^3(I)\to\C(I)$ is given by $D(f):=f'''+f'$, the null space of $D_c$ is spanned by the functions $1,\cos(t-{a}_3),\sin(t-{a}_3)$, and we have that
\Eq{*}{
  W_{{a}_1,{a}_2,{a}_3}^{1,1,1}(1,\cosh,\sinh)&=
  \begin{vmatrix}
   1 & 1 & 1 \\
   \cos({a}_1-{a}_3) & \cos({a}_2-{a}_3) & 1 \\
   \sin({a}_1-{a}_3) & \sin({a}_2-{a}_3) & 0 
  \end{vmatrix}\\
  &=\sin({a}_3-{a}_2)+\sin({a}_1-{a}_3)+\sin({a}_2-{a}_1)\\
  &=4\sin\Big(\frac{{a}_1-{a}_3}{2}\Big)\sin\Big(\frac{{a}_2-{a}_1}{2}\Big)\sin\Big(\frac{{a}_3-{a}_2}{2}\Big)\neq0.
}
Define
\Eq{5T}{
  \chi_{1,0}(t)&:=
  \frac{\cos\big(t-\frac{{a}_2+{a}_3}{2}\big)-\cos\big(\frac{{a}_2-{a}_3}{2}\big)}{\cos\big({a}_1-\frac{{a}_2+{a}_3}{2}\big)-\cos\big(\frac{{a}_2-{a}_3}{2}\big)},\\
  \chi_{2,0}(t)&:=
  \frac{\cos\big(t-\frac{{a}_1+{a}_3}{2}\big)-\cos\big(\frac{{a}_1-{a}_3}{2}\big)}{\cos\big({a}_2-\frac{{a}_1+{a}_3}{2}\big)-\cos\big(\frac{{a}_1-{a}_3}{2}\big)},\\
  \chi_{3,0}(t)&:=
  \frac{\cos\big(t-\frac{{a}_1+{a}_2}{2}\big)-\cos\big(\frac{{a}_1-{a}_2}{2}\big)}{\cos\big({a}_3-\frac{{a}_1+{a}_2}{2}\big)-\cos\big(\frac{{a}_1-{a}_2}{2}\big)}.
}
One can easily see that $\chi_{1,0},\chi_{2,0},\chi_{3,0}$ form the standard base of the kernel of $D_c$ related to the interpolational system $\big(({a}_1,{a}_2,{a}_3),(1,1,1)\big)$.
The characteristic element $\omega_c:\R\to\R$ of the null space of $D_c$ is now given as
\Eq{5To}{
\omega_c(t)=1-\cos(t) \qquad(t\in\R).
}
}

\section{The main theorem and its corollaries}

The following theorem contains the main result of this paper.

\Thm{Main}{Let $\ell,n\in\N$, $(\pmb{{a}},\pmb{n})=\big(({a}_1,\dots,{a}_l),(n_1,\dots,n_\ell)\big)\in\cH_{\ell,n}(I)$, let $c=(c_0,\dots,c_n)\in\CC^{n+1}$ with $c_n=1$ and let the differential operator $D_c:\C^n(I)\to\C(I)$ be defined by \eq{D}. Assume that the system $\big\{\chi_{\alpha,\beta}\colon\alpha\in\{1,\dots,\ell\},\beta\in\{0,\dots,n_\alpha-1\}\big\}$ is the standard base of the kernel of $D_c$ related to the interpolational system $(\pmb{{a}},\pmb{n})$.
Then, for all $f\in\C^n(I)$ and $x\in I$, we have
\Eq{Main}{
  f(x)=\sum_{\alpha=1}^\ell\sum_{\beta=0}^{n_\alpha-1} \bigg(f^{(\beta)}({a}_\alpha)
  +\int_{{a}_\alpha}^x(D_cf)(t)\cdot\omega_c^{(\beta)}({a}_\alpha-t)dt\bigg)\chi_{\alpha,\beta}(x),
}
where $\omega_c$ is the characteristic element of the kernel of $D_c$ defined by \eq{IV}.}

\begin{proof} For fixed $k\in\{0,\dots,n\}$, $\alpha\in\{1,\dots,\ell\}$, $\beta\in\{0,\dots,n_\alpha-1\}$, and $x\in I$, we can obtain
\Eq{*}{
  \int_{{a}_\alpha}^x &f^{(k)}(t)\cdot\omega_c^{(\beta)}({a}_\alpha-t)dt\\
  &=\sum_{i=1}^k\Big(f^{(k-i)}(x)\cdot\omega_c^{(\beta+i-1)}({a}_\alpha-x)
  -f^{(k-i)}({a}_\alpha)\cdot\omega_c^{(\beta+i-1)}(0)\Big)
  +\int_{{a}_\alpha}^xf(t)\cdot\omega_c^{(k+\beta)}({a}_\alpha-t)dt.
}
To check this formula, observe that the above equality is valid at $x={a}_\alpha$, and the derivatives of both sides with respect to $x\in I$ are equal everywhere. Indeed, the derivative of the right hand side equals
\Eq{*}{
  \sum_{i=1}^k\Big(f^{(k+1-i)}(x)\cdot\omega_c^{(\beta+i-1)}({a}_\alpha-x)-f^{(k-i)}(x)\cdot\omega_c^{(\beta+i)}({a}_\alpha-x)\Big)
  +f(x)\cdot\omega_c^{(k+\beta)}({a}_\alpha-x),
}
which simplifies to $f^{(k)}(x)\cdot\omega_c^{(\beta)}({a}_\alpha-x)$, which is exactly the derivative of the left hand side. Using that $\omega_c^{(\beta)}$ is in the kernel of the operator $D_c$, we can see that
\Eq{*}{
  \sum_{k=0}^nc_k\int_{{a}_\alpha}^xf(t)\cdot\omega_c^{(k+\beta)}({a}_\alpha-t)dt
  =\int_{{a}_\alpha}^xf(t)\cdot(D_c\omega_c^{(\beta)})({a}_\alpha-t)dt=0,
}
Combining the above equalities, for $\alpha\in\{1,\dots,\ell\}$, $\beta\in\{0,\dots,n_\alpha-1\}$, and $x\in I$, we get
\Eq{*}{
  f^{(\beta)}({a}_\alpha)
  &+\int_{{a}_\alpha}^x(D_cf)(t)\cdot\omega_c^{(\beta)}({a}_\alpha-t)dt\\
  &=f^{(\beta)}({a}_\alpha)
  +\sum_{k=0}^nc_k\int_{{a}_\alpha}^x f^{(k)}(t)\cdot\omega_c^{(\beta)}({a}_\alpha-t)dt \\
  &=f^{(\beta)}({a}_\alpha)
  +\sum_{k=1}^nc_k\sum_{i=1}^k\Big(f^{(k-i)}(x)\cdot\omega_c^{(\beta+i-1)}({a}_\alpha-x)
  -f^{(k-i)}({a}_\alpha)\cdot\omega_c^{(\beta+i-1)}(0)\Big)\\
  &=\sum_{k=1}^nc_k\sum_{i=1}^kf^{(k-i)}(x)\cdot\omega_c^{(\beta+i-1)}({a}_\alpha-x).
}
To see the last equality, one should use \lem{PP}. Then, for $\alpha\in\{1,\dots,\ell\}$, $\beta\in\{0,\dots,n_\alpha-1\}$, we obtain
\Eq{*}{
  \sum_{k=1}^nc_k\sum_{i=1}^k
  f^{(k-i)}({a}_\alpha)\cdot\omega_c^{(\beta+i-1)}(0)
  &=\sum_{j=0}^{n-1}f^{(j)}({a}_\alpha)\sum_{i=0}^{n-j-1}
  c_{i+j+1}\cdot\omega_c^{(\beta+i)}(0)\\
  &=\sum_{j=0}^{n-1}f^{(j)}({a}_\alpha)\cdot\delta_{j,\beta}
  =f^{(\beta)}({a}_\alpha).
}
Therefore,
\Eq{*}{
  \sum_{\alpha=1}^\ell\sum_{\beta=0}^{n_\alpha-1} &\bigg(f^{(\beta)}({a}_\alpha)
  +\int_{{a}_\alpha}^x(D_cf)(t)\cdot\omega_c^{(\beta)}({a}_\alpha-t)dt\bigg)\chi_{\alpha,\beta}(x)\\
  &=\sum_{k=1}^nc_k\sum_{i=1}^k\sum_{\alpha=1}^\ell\sum_{\beta=0}^{n_\alpha-1}f^{(k-i)}(x)\cdot\omega_c^{(\beta+i-1)}({a}_\alpha-x)\cdot\chi_{\alpha,\beta}(x).
}
For fixed $k\in\{1,\dots,n\}$, $i\in\{1,\dots,k\}$ and $y\in I$,
consider the function $F_{k,i,y}:I\to\CC$ defined by
\Eq{*}{
  F_{k,i,y}(x):=f^{(k-i)}(y)\cdot\omega_c^{(i-1)}(x-y).
}
Then, $F_{k,i,y}$ belongs to the kernel of $D_c$ (since the derivatives and translates of $\omega_c$ are included in it). On the other hand, the interpolation of $F_{k,i,y}$ with respect to the interpolational system $(\pmb{{a}},\pmb{n})$ is equal to
\Eq{*}{
  \sum_{\alpha=1}^\ell\sum_{\beta=0}^{n_\alpha-1}F_{k,i,y}^{(\beta)}({a}_\alpha)\cdot \chi_{\alpha,\beta}(x)
  =\sum_{\alpha=1}^\ell\sum_{\beta=0}^{n_\alpha-1}f^{(k-i)}(y)\cdot\omega_c^{(\beta+i-1)}({a}_\alpha-y)\cdot
  \chi_{\alpha,\beta}(x).
}
According to \lem{2}, this expression is equal to $F_{k,i,y}(x)$, i.e., for all $k\in\{1,\dots,n\}$, $i\in\{1,\dots,k\}$ and $x,y\in I$, we have the following identity
\Eq{*}{
  \sum_{\alpha=1}^\ell\sum_{\beta=0}^{n_\alpha-1}f^{(k-i)}(y)\cdot\omega_c^{(\beta+i-1)}({a}_\alpha-y)\cdot
  \chi_{\alpha,\beta}(x)=f^{(k-i)}(y)\cdot\omega_c^{(i-1)}(x-y).
}
Using this equality for $y=x$, we get that
\Eq{*}{
  \sum_{k=1}^nc_k\sum_{i=1}^k\sum_{\alpha=1}^\ell\sum_{\beta=0}^{n_\alpha-1}f^{(k-i)}(x)\cdot\omega_c^{(\beta+i-1)}({a}_\alpha-x)=\sum_{k=1}^nc_k\sum_{i=1}^kf^{(k-i)}(x)\cdot\omega_c^{(i-1)}(0)=f(x).
}
This completes the proof of the equality \eq{Main}.
\end{proof}

We now deduce the Theorem which was stated in the introduction.

\Cor{MainCl}{Let $\ell,n\in\N$, $(\pmb{{a}},\pmb{n})=\big(({a}_1,\dots,{a}_l),(n_1,\dots,n_\ell)\big)\in\cH_{\ell,n}(I)$. Then, for all $f\in\C^n(I)$ and $x\in I$, we have
\Eq{MainCl}{
  f(x)=\sum_{\alpha=1}^\ell\sum_{\beta=0}^{n_\alpha-1} \bigg(f^{(\beta)}({a}_\alpha)
  +\int_{{a}_\alpha}^x f^{(n)}(t)\cdot\frac{({a}_\alpha-t)^{n-1-\beta}}{(n-1-\beta)!}dt\bigg)H_{\alpha,\beta}(x),
}
where the polynomial $H_{\alpha,\beta}$ is defined by \eq{Hab}.}

\begin{proof} Let $c:=(0,\dots,0,1)\in\R^{n+1}$. Then the differential operator $D_c\colon\C^n(I)\to\C(I)$ and its characteristic solution $\omega_c$ are given by $D_c(f):=f^{(n)}$ and by \eq{3CE}, respectively. In this case, for $\alpha\in\{1,\dots,\ell\}$ and $\beta\in\{0,\dots,n_\alpha-1\}$, we have that $\chi_{\alpha,\beta}=H_{\alpha,\beta}$, therefore, the formula \eq{Main} of \thm{Main} reduces to \eq{MainCl}.
\end{proof}

It is not difficult to see that the equality \eq{MainCl} can equivalently be written in the form \eq{Cl} as stated in the Theorem of the introduction, therefore, \thm{Main} is indeed a generalization of it.

In the subsequent two corollaries, we describe the results related to Lagrange and Taylor interpolations, respectively. 

\Cor{3}{Let $n\in\N$, let ${a}_1,\dots,{a}_n$ pairwise distinct elements of $I$. Then, for all $f\in\C^n(I)$ and $x\in I$, we have
\Eq{*}{
  f(x)=\sum_{\alpha=1}^n\bigg(f({a}_\alpha)
  +\int_{{a}_\alpha}^x f^{(n)}(t)\frac{({a}_\alpha-t)^{n-1}}{(n-1)!}dt\bigg)
  \prod_{j\in\{1,\dots,n\}\setminus\{\alpha\}}
    \frac{x-{a}_j}{{a}_\alpha-{a}_j}.
}}

\begin{proof} We have that $\big(({a}_1,\dots,{a}_n),(1,\dots,1)\big)\in\cH_{n,n}(I)$ and let $c:=(0,\dots,0,1)\in\R^{n+1}$. Then the differential operator $D_c\colon\C^n(I)\to\C(I)$ is given by $D_c(f):=f^{(n)}$. For $i\in\{1,\dots,n\}$, define $\chi_{i,0}:\R\to\R$ by \eq{3SB}. According to \exa{3}, these polynomials form the standard base of the null space of $D_c$ related to the interpolational system $\big(({a}_1,\dots,{a}_n),(1,\dots,1)\big)$. The characteristic element $\omega_c:\R\to\R$ of $\ker(D_c)$ now is given by \eq{3CE}. With these expressions, the equality \eq{Main} simplifies to the asserted equality.
\end{proof}

In the next result we deduce the generalized Taylor theorem obtained in the paper \cite{AliPal22}.

\Cor{4}{Let $n\in\N$, ${a}\in I$, let $c=(c_0,\dots,c_n)\in\CC^{n+1}$ with $c_n:=1$. Define the differential operator $D_c:\C^n(I)\to\C(I)$ by \eq{D} and let $\omega_c:\R\to\R$ denote the characteristic element of $\ker(D_c)$. Then, for all $f\in\C^n(I)$ and $x\in I$, 
\Eq{*}{
  f(x)=\sum_{\beta=0}^{n-1} f^{(\beta)}({a})\bigg(\sum_{i=0}^{n-1-\beta}c_{i+\beta+1}\omega_c^{(i)}(x-{a})\bigg)
  +\int_{{a}}^x(D_cf)(t)\omega_c(x-t)dt.
}}

\begin{proof}
According to \exa{4+}, the functions $\chi_{1,0},\dots,\chi_{1,n-1}$ defined by \eq{4SB+} the standard base of the kernel of $D_c$ related to the interpolational system $({a},n)$. Therefore, according to \thm{Main}, we have that
\Eq{TF4}{
  f(x)&=\sum_{\beta=0}^{n-1} \bigg(f^{(\beta)}({a})
  +\int_{{a}}^x(D_cf)(t)\cdot\omega_c^{(\beta)}({a}-t)ds\bigg)\sum_{i=0}^{n-1-\beta}c_{i+\beta+1}\omega_c^{(i)}(x-{a})\\
  &=\sum_{\beta=0}^{n-1} f^{(\beta)}({a})\sum_{i=0}^{n-1-\beta}c_{i+\beta+1}\omega_c^{(i)}(x-{a})
  +\int_{{a}}^x(D_cf)(t)\sum_{\beta=0}^{n-1}\sum_{i=0}^{n-1-\beta}c_{i+\beta+1}\omega_c^{(i)}(x-{a})\omega_c^{(\beta)}({a}-t)dt.
}
To complete the proof of this theorem, using \thm{AF} with $\omega:=\omega_c$ and $u:=x-{a}$, $v:={a}-t$, for all $x,t\in\R$, we can get that
\Eq{*}{
  \sum_{\beta=0}^{n-1}\sum_{i=0}^{n-1-\beta}c_{i+\beta+1}\omega_c^{(i)}(x-{a})\omega_c^{(\beta)}({a}-t)=\omega_c(x-t).
}
In view of this equality, the previous equality simplifies to the asserted formula \eq{TF4}.
\end{proof}

The following result is the classical Taylor theorem with integral remainder term. It directly follows from the previous corollary.

\Cor{TF4}{Let $n\in\N$ and let ${a}\in I$. Then, for all $f\in\C^n(I)$ and $x\in I$, we have
\Eq{*}{
  f(x)=\sum_{\beta=0}^{n-1}f^{(\beta)}({a})\frac{(x-{a})^\beta}{\beta!}+\int_{{a}}^x f^{(n)}(t)\frac{(x-t)^{n-1}}{(n-1)!}dt.
}}

\begin{proof} Let $c:=(0,\dots,0,1)\in\R^{n+1}$. Then the differential operator $D_c\colon\C^n(I)\to\C(I)$ is given by $D(f):=f^{(n)}$. The characteristic element $\omega_c:\R\to\R$ of $\ker(D_c)$ is given by \eq{3CE}. Applying the formula of \cor{4}, the result follows.
\end{proof}

\Cor{1H}{Let ${a}_1,{a}_2\in I$ such that ${a}_1\neq{a}_2$.
Then, for all $f\in\C^2(I)$ and $x\in I$, we have
\Eq{*}{
  f(x)&=\bigg(f({a}_1)+\int_{{a}_1}^x(f''-f)(t)\cdot\sinh({a}_1-t)dt\bigg)\frac{\sinh(x-{a}_2)}{\sinh({a}_1-{a}_2)}\\
  &\quad+\bigg(f({a}_2)+\int_{{a}_2}^x(f''-f)(t)\cdot\sinh({a}_2-t)dt\bigg)\frac{\sinh(x-{a}_1)}{\sinh({a}_2-{a}_1)}.
}}

\begin{proof} Let us apply \thm{Main} to the interpolation system $(({a}_1,{a}_2),(1,1))\in\cH_{2,2}(I)$ and to the vector $c=(-1,0,1)$. Then the differential operator $D_c$ is given as $D_c(f)=f''-f$ and, according to \exa{1H},  the standard base $\{\chi_{1,0},\chi_{2,0}\}$ for the kernel of $D_c$ is given by \eq{1H}. One can also check that the characteristic solution $\omega_c$ is now given by \eq{1Ho}. With these expressions, the equality \eq{Main} simplifies to the asserted equality.
\end{proof}

\Cor{1T}{Let ${a}_1,{a}_2\in I$ such that $\sin({a}_1-{a}_2)\neq0$. Then, for all $f\in\C^2(I)$ and $x\in I$, we have
\Eq{*}{
  f(x)&=\bigg(f({a}_1)+\int_{{a}_1}^x(f''-f)(t)\cdot\sin({a}_1-t)dt\bigg)\frac{\sin(x-{a}_2)}{\sin({a}_1-{a}_2)}\\
  &\quad+\bigg(f({a}_2)+\int_{{a}_2}^x(f''-f)(t)\cdot\sin({a}_2-t)dt\bigg)\frac{\sin(x-{a}_1)}{\sin({a}_2-{a}_1)}.
}}

\begin{proof} We apply \thm{Main} to the interpolation system $(({a}_1,{a}_2),(1,1))\in\cH_{2,2}(I)$, to the vector $c=(1,0,1)$, and we use the formulae obtained in \exa{1T}.
Then the assertion follows analogously as in the proof of \cor{1H}. 
\end{proof}

\Cor{2HT}{Let ${a}_1,{a}_2\in I$ and assume that $\cosh({a}_1-{a}_2)\cos({a}_1-{a}_2)\neq 1$.
Denote $K:=2-2\cosh({a}_1-{a}_2)\cos({a}_1-{a}_2)$. Then, for all $x\in I$,
\Eq{*}{
  f(x)&=\frac{1}{K}\bigg[\bigg(f({a}_1)+\frac12\int_{{a}_1}^x (f^{(4)}-f)(t)\big(\sinh({a}_1-t)-\sin({a}_1-t)\big)dt\bigg)\\
  &\qquad\qquad\times\begin{vmatrix}
   \cosh({a}_1-{a}_2)-\cos({a}_1-{a}_2) 
   & \sinh(x-{a}_2)-\sin(x-{a}_2) \\
   \sinh({a}_1-{a}_2)+\sin({a}_1-{a}_2) 
   & \cosh(x-{a}_2)-\cos(x-{a}_2)
   \end{vmatrix}\\
   &\qquad+\bigg(f'({a}_1)+\frac12\int_{{a}_1}^x (f^{(4)}-f)(t)\big(\cosh({a}_1-t)-\cos({a}_1-t)\big)dt\bigg)\\
  &\qquad\qquad\times\begin{vmatrix}
   \cosh({a}_1-{a}_2)-\cos({a}_1-{a}_2) 
   & \cosh(x-{a}_2)-\cos(x-{a}_2) \\
   \sinh({a}_1-{a}_2)-\sin({a}_1-{a}_2) 
   & \sinh(x-{a}_2)-\sin(x-{a}_2)
   \end{vmatrix}\\
   &\qquad+\bigg(f({a}_2)+\frac12\int_{{a}_2}^x (f^{(4)}-f)(t)\big(\sinh({a}_2-t)-\sin({a}_2-t)\big)dt\bigg)\\
   &\qquad\qquad\times\begin{vmatrix}
   \cosh({a}_1-{a}_2)-\cos({a}_1-{a}_2) 
   & -\sinh(x-{a}_1)+\sin(x-{a}_1) \\
   \sinh({a}_1-{a}_2)+\sin({a}_1-{a}_2) 
   & \cosh(x-{a}_1)-\cos(x-{a}_1)
   \end{vmatrix}\\
   &\qquad+\bigg(f'({a}_2)+\frac12\int_{{a}_2}^x (f^{(4)}-f)(t)\big(\cosh({a}_2-t)-\cos({a}_2-t)\big)dt\bigg)\\
  &\qquad\qquad\times\begin{vmatrix}
   \cosh({a}_1-{a}_2)-\cos({a}_1-{a}_2) 
   & -\cosh(x-{a}_1)+\cos(x-{a}_1) \\
   \sinh({a}_1-{a}_2)-\sin({a}_1-{a}_2) 
   & \sinh(x-{a}_1)-\sin(x-{a}_1)
   \end{vmatrix}\bigg].
}
}

\begin{proof} We apply \thm{Main} to the interpolational system $\big(({a}_1,{a}_2),(2,2)\big)\in\cH_{2,4}(I)$ and let $c:=(-1,0,0,0,1)$. Then the $4$th-order differential operator $D_c\colon\C^4(I)\to\C(I)$ is defined by $D_c(f):=f^{(4)}-f$. Define the functions $\{\chi_{1,0},\chi_{1,1},\chi_{2,0},\chi_{2,1}\}$ by \eq{2HT}. Then, according to \exa{2HT}, $\{\chi_{1,0},\chi_{1,1},\chi_{2,0},\chi_{2,1}\}$ is the standard base of the kernel of $D_c$ related to the interpolational system $\big(({a}_1,{a}_2),(2,2)\big)$. The characteristic element $\omega_c:\R\to\R$ of $\ker(D_c)$ is given by \eq{2HTo}. With these expressions, the equality \eq{Main} simplifies to the asserted equality.
\end{proof}

\Cor{5H}{Let ${a}_1,{a}_2,{a}_3\in I$ be pairwise distinct point. Then, for all $f\in\C^3(I)$ and $x\in I$, we have
\Eq{*}{
  f(x)&=\bigg(f({a}_1)+\int_{{a}_1}^x(f'''-f')(t)\cdot(\cosh({a}_1-t)-1)dt\bigg)\frac{\cosh\big(x-\frac{{a}_2+{a}_3}{2}\big)-\cosh\big(\frac{{a}_2-{a}_3}{2}\big)}{\cosh\big({a}_1-\frac{{a}_2+{a}_3}{2}\big)-\cosh\big(\frac{{a}_2-{a}_3}{2}\big)}\\
  &\quad+\bigg(f({a}_2)+\int_{{a}_2}^x(f'''-f')(t)\cdot(\cosh({a}_2-t)-1)dt\bigg)\frac{\cosh\big(x-\frac{{a}_1+{a}_3}{2}\big)-\cosh\big(\frac{{a}_1-{a}_3}{2}\big)}{\cosh\big({a}_2-\frac{{a}_1+{a}_3}{2}\big)-\cosh\big(\frac{{a}_1-{a}_3}{2}\big)}\\
  &\quad+\bigg(f({a}_3)+\int_{{a}_3}^x(f'''-f')(t)\cdot(\cosh({a}_3-t)-1)dt\bigg)\frac{\cosh\big(x-\frac{{a}_1+{a}_2}{2}\big)-\cosh\big(\frac{{a}_1-{a}_2}{2}\big)}{\cosh\big({a}_3-\frac{{a}_1+{a}_2}{2}\big)-\cosh\big(\frac{{a}_1-{a}_2}{2}\big)}.
}}

\begin{proof} We apply \thm{Main} to the interpolational system $\big(({a}_1,{a}_2,{a}_3),(1,1,1)\big)\in\cH_{3,3}(I)$ and let 
$c:=(0,-1,0,1)$. Then the $3$rd-order differential operator $D_c\colon\C^3(I)\to\C(I)$ is given as $D(f):=f'''-f'$. According to \exa{5H}, the functions $\chi_{1,0},\chi_{2,0},\chi_{3,0}$ defined in \eq{5H} form the standard base of the kernel of $D_c$ related to the interpolational system $\big(({a}_1,{a}_2,{a}_3),(1,1,1)\big)$. The characteristic element $\omega_c:\R\to\R$ of the null space of $D_c$ is now given by the equality \eq{5Ho}, therefore, the formula \eq{Main} yields the asserted equality.
\end{proof}

\Cor{5T}{Let ${a}_1,{a}_2,{a}_3\in I$ with $
\sin\big(\frac{{a}_1-{a}_3}{2}\big)\sin\big(\frac{{a}_2-{a}_1}{2}\big)\sin\big(\frac{{a}_3-{a}_2}{2}\big)\neq0$. Then, for all $f\in\C^3(I)$ and $x\in I$, we have
\Eq{*}{
  f(x)&=\bigg(f({a}_1)+\int_{{a}_1}^x(f'''+f')(t)\cdot(1-\cos({a}_1-t))dt\bigg)\frac{\cos\big(x-\frac{{a}_2+{a}_3}{2}\big)-\cos\big(\frac{{a}_2-{a}_3}{2}\big)}{\cos\big({a}_1-\frac{{a}_2+{a}_3}{2}\big)-\cos\big(\frac{{a}_2-{a}_3}{2}\big)}\\
  &\quad+\bigg(f({a}_2)+\int_{{a}_2}^x(f'''+f')(t)\cdot(1-\cos({a}_2-t))dt\bigg)\frac{\cos\big(x-\frac{{a}_1+{a}_3}{2}\big)-\cos\big(\frac{{a}_1-{a}_3}{2}\big)}{\cos\big({a}_2-\frac{{a}_1+{a}_3}{2}\big)-\cos\big(\frac{{a}_1-{a}_3}{2}\big)}\\
  &\quad+\bigg(f({a}_3)+\int_{{a}_3}^x(f'''+f')(t)\cdot(1-\cos({a}_3-t))dt\bigg)\frac{\cos\big(x-\frac{{a}_1+{a}_2}{2}\big)-\cos\big(\frac{{a}_1-{a}_2}{2}\big)}{\cos\big({a}_3-\frac{{a}_1+{a}_2}{2}\big)-\cos\big(\frac{{a}_1-{a}_2}{2}\big)}.
}}

\begin{proof} The proof of this corollary is completely analogous to that of \cor{5H}, therefore it is omitted.
\end{proof}


\end{document}